\documentclass[12pt]{amsart}
\usepackage{amsfonts, amsmath, amsthm}
\usepackage{graphicx}

\newtheorem{pro}{Proposition}[section]
\newtheorem{thm}[pro]{Theorem}
\newtheorem{lem}[pro]{Lemma}

\newtheorem{cnj}[pro]{Conjecture}

\newtheorem{cor}[pro]{Corollary}

\theoremstyle{definition}
\newtheorem{dfn}[pro]{Definition}
\newtheorem{dfns}[pro]{Definitions}

\newtheorem{rmkk}[pro]{Remark}
\newtheorem{rmkks}[pro]{Remarks}

\theoremstyle{remark}

\newcommand{\s}{\Sigma}

\newcommand{\valchi}{4 - 2(m+c)g}
\newcommand{\valchin}{4 - 2gn}
\newcommand{\val}{\max\{2gn-2,2(h+n-1)\}}
\newcommand{\kk}{\mathcal{K}_{g,n}}

\title[Knot exteriors with additive Heegaard genus]
{Knot exteriors with additive Heegaard genus and Morimoto's Conjecture}

\date{\today} \address{Department of Mathematics, Nara Women's University
Kitauoya Nishimachi, Nara 630-8506, Japan} \address{Department of mathematical
Sciences, University of Arkansas, Fayetteville, AR 72701}
\email{tsuyoshi@cc.nara-wu.ac.jp} \email{yoav@uark.edu} \author{Tsuyoshi
Kobayashi} \author{Yo'av Rieck}
\begin{document}


\date{\today}%
%\dedicatory{}
%\commby{}%
% ----------------------------------------------------------------
\begin{abstract}
Given integers $\{g_i \geq 2 \}_{i=1}^n$ we prove that there exists infinitely
may knots $K_i \subset S^3$ so that $g(E(K_i)) = g_i$ and $g(E(\#_{i=1}^n
  K_i) = \s_{i=1}^n g(E(K_i))$.  (Here, $E(\cdot)$ denotes the exterior and
  $g(\cdot)$ the Heegaard genus.)
Together with \cite[Theorem~1.5]{crelle}, this proves the existence of
counterexamples to Morimoto's Conjecture \cite{mc}.
\end{abstract}
\maketitle
% ------------------------------------------------------

\section{Introduction and statements of results}
\label{intro}

Let $K_i$ ($i=1,2$) be knots in the 3-sphere $S^3$, and let $K_1 \# K_2$ be
their connected sum.  We use the notation $t(\cdot)$, $E(\cdot)$, and
$g(\cdot)$ to denote tunnel number, exterior, and Heegaard genus respectively
(we follow the definitions and notations given in \cite{cag}).  It is well
known that the union of a tunnel system for $K_1$, a tunnel system for $K_2$,
and a tunnel on a decomposing annulus for $K_1 \# K_2$ forms a tunnel system
for $K_1 \# K_2$.  Therefore:
$$t(K_1 \# K_2) \leq t(K_1) + t(K_2) +1.$$
Since (for any knot $K$) $t(K) = g(E(K)) - 1$, this gives:
\begin{equation}
\label{eq:upper-bound}  
g(E(K_1 \# K_2)) \leq g(E(K_1)) + g(E(K_2)).
\end{equation}

We say that a knot $K$ in a closed orientable manifold $M$ admits a $(g,n)$
position if there exists a genus $g$  Heegaard surface $\s \subset M$,
separating $M$ into the handlebodies $H_1$ and $H_2$, so that $H_i \cap K$
($i=1,2$) consists of $n$ arcs that are simultaneously parallel into $\partial
H_i$.  We say that $K$ admits a $(g,0)$ position if $g(E(K)) \leq g$.  
Note that if $K$ admits a $(g,n)$ position then $K$ admits both a $(g,n+1)$
position and a $(g+1,n)$ position.  

\begin{rmkk}
\label{rmk:g,0}
The definition given in \cite{jt} for $(g,n)$ position with $n \geq 1$ is
identical to our definition.  However, in \cite{jt} $K$ is said to admit a
$(g,0)$ position if $K$ is isotopic into a  genus $g$ Heegaard surface for $M$.
Thus, if $K$ admit a $(g,0)$ position in the sense of \cite{jt} and $g(X) > g$ 
then $K$ admits a $(g,1)$ position in our sense. For example, a non-trivial
torus knot in $S^3$ is called $(1,0)$ in \cite{jt} and $(1,1)$ here.  ({\it
  Cf.} \cite[Remark~2.4]{mms}.)
\end{rmkk}

It is known \cite[Proposition~1.3]{mc} that if $K_i$ ($i=1$ or 2) admits a
$(t(K_i),1)$ position then equality does not hold in
Inequality~(\ref{eq:upper-bound}):
\begin{equation}
\label{eq:upper-bound-inequality}  
g(E(K_1 \# K_2)) < g(E(K_1)) + g(E(K_2)).
\end{equation}
Morimoto proved that if $K_1$ and $K_2$ are m-small knots in $S^3$ then the converse
holds \cite[Theorem~1.6]{mc}.  
This result was generalized to arbitrarily many m-small knots in general
manifolds by the authors \cite{cag}.  Morimoto conjectured that the converse
holds in general \cite[Conjecture~1.5]{mc}: 
\begin{cnj}[Morimoto's Conjecture]
\label{cnj:mc}
Given knots $K_1,\ K_2 \subset S^3$, $g(E(K_1 \# K_2)) < g(E(K_1)) +
g(E(K_2))$ if and only if $K_i$ admits a $(t(K_i),1)$ position (for $i=1$ or
$i=2$).
\end{cnj}

\begin{rmkk}

We note that Morimoto stated the above conjecture in terms of 1-bridge genus 
$g_1(K)$. It is easy to see that the Conjecture~1.5 of \cite{mc} is equivalent
to the statement above.
\end{rmkk}

In \cite{crelle} the authors showed that certain conditions imply existence of
counterexamples to Morimoto's Conjecture.  One such condition is the existence
of an m-small knot $K$ that does not admit a $(t(K),2)$ position.  We asked
\cite[Question~1.9]{crelle} if there exists a knot $K$ with $g(E(K)) = 2$ that
does not admit a $(1,2)$ position; this question was answered affirmatively by
Johnson and  Thompson \cite[Corollary~2]{jt}, who showed that for any $n$
there exist infinitely many knots with $g(E(K))=2$ not admitting a $(1,n)$
position.  At about the same time Minsky, Moriah and Schleimer
\cite[Theorem~4.2]{mms} proved a more general result, showing that for any
integers $g \geq 2$, $n \geq 1$ there exist infinitely many knots with
$g(E(K)) = g$ that do not admit a $(g-1,n)$ position (more precisely, this
follows from \cite[Theorem~3.1]{mms} and Proposition~\ref{pro:jt}  
below).  Although it is not known if any of these examples are m-small, in
this paper we show that some of these examples have the property described in
the theorem below, that also implies existence of counterexamples to
Morimoto's Conjecture.

\begin{thm}
\label{thm:additive}
Given integers $g \geq 2$ and $n \geq 1$, there exists a family of knots in
$S^3$ (denoted $\kk$) with the following properties:
\begin{enumerate}
\item For each $h$ with $2 \leq h \leq g$, there exists infinitely many knots 
$K \in \kk$ with $g(E(K)) = h$.
\item For any collection of knots $K_1,\dots,K_m \in \kk$ 
(possibly, $K_i = K_j$ for $i \ne j$) with $m \leq n$,
  $$g(E(\#_{i=1}^m K_i)) = \s_{i=1}^m g(E(K_i)).$$
\end{enumerate}
Moreover, for each $g$, we have:
$$\cap_{n=1}^\infty \kk = \emptyset.$$
\end{thm}

\begin{rmkks}
\label{rmk}
  \begin{enumerate}
  \item The knots in $\kk$ need not be prime.  In fact, it is clear from the
  definition of $\kk$ that if $K \in \mathcal{K}_{g,pn}$ then $pK \in
  \mathcal{K}_{pg,n}$ ($pK$ is defined in Definitions~\ref{defns}).  We do not  
  know if   $\kk$ contains a knot of the form $pK$ (for $p >1$) when $g$ is
  prime.  
%   \item A $(t(K),1)$ position is used extensively in the literature.  One of
%    the the contributions of this paper is using the value of $b$ for
%    $(t(K),b)$ position when $b > 1$.  This is refined further in$rieck
%    \cite{strongHH} (at the time of writing in preparation).
  \item Existence of knots $K_1,\ K_2$ with $g(E(K_1 \# K_2)) = g(E(K_1)) +
        g(E(K_2))$ is known from \cite{moriah-rubinstein} and \cite{msy}.
        Theorem~\ref{thm:additive} is new in the following ways:
            \begin{enumerate}
                  \item It is the first time that the connected sum of 
                  more than two knots are shown to
                  have additive Heegaard genus.
                  \item The proof in \cite{moriah-rubinstein} uses minimal
                  surfaces in hyperbolic manifolds and in \cite{msy} quantum
                  invariants.  Our proof is purely topological.
  \end{enumerate}
\item The sets $\kk$ are not uniquely defined; for example, we can remove any
  finite set from $\kk$.  However, for any sets $\kk$ fulfilling
  Theorem~\ref{thm:additive}~(1) and (2), we have that $\cap_n \kk =
  \emptyset$.  
\end{enumerate}
\end{rmkks}

A knot $K \subset M$ is called {\it admissible} (see \cite{crelle}) if 
$g(E(K)) > g(M)$.  Thus any knot $K \subset S^3$ is admissible.  By
\cite[Theorem~1.2]{crelle} for any admissible knot $K$ there exists $N$ so
that if $n \geq N$ then $g(E(nK)) < n g(E(K))$.   In contrast to that we have:

\begin{cor}
  \label{cor:additive}
Given integers $g \geq 2$ and $n \geq 1$, there exist an infinitely many knots
$K \subset S^3$ so that $g(E(K))=g$ and for any $m \leq n$, $g(E(mK)) = mg$.
\end{cor}

\begin{proof}
For $K \in \kk$ with $g(E(K)) = g$ we have $g(E(nK)) = ng$.
\end{proof}

\begin{rmkk}
By \cite[Proposition~1.7]{crelle}, a knot $K$ with $g(E(K))=g$ and $g(E(nK)) =
ng$ cannot admit a $(g(X) - 1, n-1)$ position.     
\end{rmkk}

Another consequence of Corollary~\ref{cor:additive} is:

\begin{cor}
\label{cor:mc}
There exists a counterexample to Morimoto's Conjecture, specifically, there
exist knots $K_1, \ K_2 \subset S^3$ so that $K_i$ does not admit a
$(t(K_i),1)$ position ($i=1,2$), and (for some integer $m$) $g(E(K_1)) = 4$,
$g(E(K_2)) = 2(m-2)$, and $g(E(K_1 \# K_2)) < 2m$.
\end{cor}

\begin{proof}[Proof of Corollary~\ref{cor:mc}]
This argument was originally given in \cite[Theorem~1.4]{crelle}.   We outline
it here for completeness.  Let $K$ be a knot as in
Corollary~\ref{cor:additive}, for $g=2$ and $n=3$.   By
\cite[Theorem~1.2]{crelle}, for some $m > 1$, $g(E(mK)) < m g(E(K)) = 2m$.
Let $m$ be the minimal number with that property.  By
Corollary~\ref{cor:additive}, $m \geq 4$.    Hence $g(E(2K)) = 2g(E(K)) = 4$.
By the minimality of $m$, $g(E((m-2)K)) = (m-2)g(E(K)) = 2(m-2)$.

Let $K_1 = 2K$ and $K_2 = (m-2)K$.  Note that $K_1 \# K_2 = mK$.
We have seen: 
\begin{enumerate}
\item $g(E(K_1)) = 4$.
\item $g(E(K_2)) = 2(m-2)$.
\item $g(E(K_1 \# K_2)) < 2m$.
\end{enumerate}

We claim that $K_1$ does not admit a $(t(K_1),1)$ position; assume for a
contradiction it does.  By  Inequality~(\ref{eq:upper-bound-inequality})
and the above (1),  we would have that $g(E(3K)) = g(E(K_1 \# K)) < g(E(K_1))
+ g(E(K)) = 6$, contradicting our choice of $K$.

We claim that $K_2$ does not admit a $(t(K_2),1)$ position; assume for a
contradiction it does.   Then by
Inequality~(\ref{eq:upper-bound-inequality})  and the above (2),
$g(E((m-1)K)) < g(E((m-2)K)) + g(E(K)) = (m-1)g(E(K))$, contradicting
minimality of $m$.
\end{proof}

We note that $K_1$ and $K_2$ are composite knots.  This leads Moriah
\cite[Conjecture~7.14]{moriah-s} to conjecture that if $K_1$ and $K_2$ are
prime then Conjecture~\ref{cnj:mc} holds.

\medskip

\noindent{\bf Outline.} Section~\ref{sec:props} is devoted to three
propositions necessary for the  proof of Theorem~\ref{thm:additive}:
Proposition~\ref{pro:bridge} that relates  strongly irreducible Heegaard
splittings and bridge position, Proposition~\ref{pro:essential} that relates
Essential surfaces and the {\it distance} of Heegaard splitting
(Proposition~\ref{pro:essential} is exactly Theorem~3.1 of
\cite{scharlemann-2004}), and Proposition~\ref{pro:jt} which relates bridge
position and distance of Heegaard splittings (Proposition~\ref{pro:jt} is
based on and extends  Theorem~1 of \cite{jt}).   In Section~\ref{sec:genus} we
calculate the genera of certain manifolds that we denote by $X(m)^{(c)}$.  In
Section~\ref{sec:proof} we prove Theorem~\ref{thm:additive}.

\begin{rmkk}
The reader may wish to read \cite{small}, where an easy argument is given for
a special case of Corollary~\ref{cor:additive}, namely, $g=2$ and $n=3$.
Note that this special case is sufficient for Corollary~\ref{cor:mc};
\cite{small} can be used as an introduction to the ideas in the current
paper. 
\end{rmkk}

\section{Decomposing $X^{(c)}$.}
\label{sec:props}

In this and the following sections, we adopt the following notations.

\begin{dfns}
\label{defns}
Let $K$ be a knot in a closed orientable manifold and $X$ its exterior.  Let
$n \geq 1$ be an integer.
  \begin{enumerate}
  \item The connected sum of $n$ copies of $K$ is denoted by $nK$ and its
  exterior by $X(n)$.
  \item For an integers $c \geq 0$ and $n \geq 1$ we denote by $X(n)^{(c)}$
    the manifold obtained by drilling $c$ curves out of 
    $X(n)$ that are simultaneously parallel to meridians of $nK$.  For
    convenience, we denote $X(1)^{(c)}$ by $X^{(c)}$.  
(Note that $X^{(0)}=X$, and $X(n)^{(0)}=X(n)$.)
  \end{enumerate}
\end{dfns}

\begin{pro}
\label{pro:bridge}
Let $X$, $X^{(c)}$ be as above and $g \geq 0$ an integer.  
Suppose $c>0$, and $X^{(c)}$
admits a strongly irreducible Heegaard surface of genus $g$. Then one of the
following holds:
\begin{enumerate}
\item $X$ admits an essential surface $S$ with $\chi(S) \geq 4 - 2g$.
\item For some $b$, $c \leq b \leq g$, $K$ admits a $(g-b,b)$ position.
\end{enumerate}
\end{pro}

\begin{rmkks}
  \begin{enumerate}
  \item If $c > g$ then conclusion~(1) holds.
  \item   Compare the proof to \cite[Theorem~3.8]{schsch}.
  \end{enumerate}
\end{rmkks}

\begin{proof}[Proof of Proposition~\ref{pro:bridge}]
Let $C_1 \cup_\s C_2$ be a genus $g$ strongly irreducible Heegaard splitting of
$X^{(c)}$.

Since $c>0$, $X^{(c)}$ admits an essential torus $T$ that gives the decomposition
$X^{(c)} = X' \cup_T Q^{(c)}$, where $X' \cong X$ and $Q^{(c)}$ is a $c$-times
punctured annulus cross $S^1$.  Since $T$ is incompressible and $\s$ is strongly
irreducible, we may isotope $\s$ so that every component of $\s \cap T$ is
essential in both surfaces.  Isotope $\s$ to minimize $|\s \cap T|$ subject to
that constraint.  Denote $\s \cap X'$ by $\s_X$ and $\s \cap Q^{(c)}$ by $\s_Q$.
Note that (by essentiality of $T$) $\s \cap T \neq \emptyset$ and (by
minimality) no component of $\s_X$ (resp. $\s_Q$) is boundary parallel in $X'$
(resp. $Q^{(c)}$).  By the argument of \cite[Claim~4.5]{cag} we
may assume that $\s_X$ is connected and compresses into both sides in $X'$ and
$\s_Q$ is incompressible in $Q^{(c)}$, for otherwise Conclusion~(1) holds.

Every component of $\s_Q$ is a vertical annulus (see, for example,
\cite[VI.34]{jaco}).  Hence, $\partial \s_X$ consists of meridians of $K$.
For $i=1,2$, let $\s_i$ be the surface obtained by simultaneously compressing
$\s_X$ maximally into $C_i \cap X$.  Then the argument of Claim~6 (page 248)
of \cite{KRlocal}  shows that every component of $\s_i$ is incompressible.
Hence, we may assume that every component is  a boundary parallel annulus in
$X'$ or a 2-sphere (for otherwise Conclusion~(1) holds).  Denote that number
of annuli by $b$ (note that $b = \frac{1}{2}|\partial \s_X|$ and is the same
for $\s_1$ and $\s_2$).  Denote the solid tori that define that boundary
parallelism of $\s_i$ by $N_{i,1},\dots,N_{i,b}$.

\medskip

\noindent{Claim.}  For $i=1,2$, $N_{i,1},\dots,N_{i,b}$ are mutually
disjoint.

\begin{proof}[Proof of claim.]
Assume for a contradiction that two components (say $N_{i,1}$ and $N_{i,2}$)
intersect, say $N_{i,2} \subset N_{i,1}$.  By construction $\s_X$ is a
connected surface obtained by tubing the annuli $\s_i$ and (possibly empty)
collection of 2-spheres into one side, therefore the tubes are all contained
in $\mbox{cl}(N_{i,1} \setminus N_{i,2})$, and we see that (for all $j \ge 2$)
$N_{i,j} \subset N_{i,1}$.  This shows that $\s$ is isotopic into $Q^{(c)}$,
hence $T \subset C_1$ or $T \subset C_2$.  Since $T$ is essential, this is
impossible.  This proves the claim.
\end{proof}

By the claim, $C_i \cap X'$ is obtained from $N_{i,1},\dots, N_{i,b}$ and a
(possibly empty) collection of 3-balls by attaching 1-handles.  This implies that
$C_i \cap X'$ is obtained from a handlebody  $H$ (say of genus $h$) by
removing a regular neighborhood of $b$ trivial arcs, say
$\gamma_{i,1}\dots,\gamma_{i,b}$, where $N_{i,j} \cap T$ corresponds to the
frontier of the regular neighborhood of $\gamma_{i,j}$ ($j=1,\dots,b$).  Since
every component of $\s_Q$ is an annulus, $\chi(\s_X) = \chi(\s)$.  $\partial
H$ is obtained by capping $\s_X$ off with $2b$ disks, hence $\chi(\partial H)
= \chi(\s) + 2b$; this shows that $h = g - b$.

We obtained a $(g-b,b)$ position for $K$, and to complete the proof we need to
show that $b \geq c$.  Suppose for a contradiction that $b<c$.  Note that
$\s_Q$ consists of $b$ vertical annuli.  Since $\partial \s_Q \subset T$, we
see that  $\s_Q$ separates $Q^{(c)}$ into $b+1$ components.  Note that
$\partial X^{(c)}$ consists of $c+1$ tori; thus if $b < c$ then two components
of $\partial Q^{(c)}$ are in the same component of $Q^{(c)}$ cut open along
$\s_Q$.  It is easy to see that there is a vertical annulus connecting  these
tori which is disjoint from $\s$.  Hence this annulus is contained in the
compression body $C_i$ and connects components of $\partial_- C_i$ for $i=1$
or $2$, a contradiction (for the notation $\partial_- C_i$,  see, for example,
\cite{cag}).   This contradiction completes the proof of
Proposition~\ref{pro:bridge}.
\end{proof}

\begin{dfn}[Hempel \cite{hempel}]
Let $H_1 \cup_\s H_2$ be a Heegaard splitting. The {\it distance} of $\s$,
denoted $d(\s)$, is the least integer
$d$ so that there exist meridian disks $D_i \subset H_i$ ($i=1,2$) and
essential curves $\gamma_0,\dots,\gamma_d \subset \s$ so that $\gamma_0 =
\partial D_1$, $\gamma_d = D_2$, and $\gamma_{i-1} \cap \gamma_i = \emptyset$
($i=1,\dots,d$).  If $\s$ is the trivial Heegaard splitting of a compression
body (that is, $\s$ is boundary parallel) this definition does not apply,
since on one side of $\s$ there are no meridional disks.  In that case, we
define $d(\s)$ to be zero.  
\end{dfn}

The properties of knots with exteriors of high distance that we need are given
in the next two propositions:

\begin{pro}
\label{pro:essential}
Let $K$ be a knot and $d \geq 0$ an integer. Suppose $X$ admits a Heegaard
splitting with distance greater than $d$.  Then $X$ does not admit a connected
essential surface $S$ with $\chi(S) \geq 2-d $.
\end{pro}

\begin{proof}
This is Theorem~3.1 of \cite{scharlemann-2004}.
\end{proof}

Our next proposition is a combination of Theorem~1 of \cite{jt} and
Corollary~3.5 of \cite{schtom}:

\begin{pro}
\label{pro:jt}  
Let $K \subset M$ be a knot and $p, \ q$ integers so that $K$ admits a $(p,q)$
position. 

If $p < g(X)$ then any Heegaard splitting for $X$ has distance at most
$2(p+q)$.
\end{pro}

\begin{proof}
Recall from Remark~\ref{rmk:g,0} that our definition of $(p,q)$ position is
not quite the same as \cite{jt}.  As explained in Remark~\ref{rmk:g,0} either
$K$ admits a $(p,q)$ position in the sense of \cite{jt} or $q=1$ and $K$
admits a $(p,0)$ position in the sense of \cite{jt}.  In the former case,
Proposition~\ref{pro:jt} is exactly Theorem~1 of \cite{jt}.  Thus we may
assume:
\begin{enumerate}
\item $q=1$.
\item $K$ admits a $(p,0)$ position in the sense of \cite{jt}, that is, $M$
  admits a Heegaard splitting of genus $p$ (say $H_1 \cup_\s H_2$) so that $K$
  is isotopic into $\s$.
\end{enumerate}

We base our analysis on
\cite{rieck}\cite{rieck-sedgwick1}\cite{rieck-sedgwick2}.  After isotoping $K$
into $\s$, let $N = N_{H_1}(K)$ be a neighborhood of $K$ in $H_1$.  Then $N$
is a solid torus and $K \subset \partial N$ a longitude.   Let $\Delta$ be a
meridian disk of $N$ that intersects $K$ in one point.   Let $\alpha \subset
\Delta$ be a properly embedded arc with  $\partial \alpha \subset (\Delta \cap
\s)$, so that $K \cap (\Delta \cap \s)$ separates the points of  $\partial
\alpha$.  Let $K'$ be a copy of $K$ pushed slightly into $H_1$, so that  $K'
\cap \Delta$ is a single point and $\alpha$ does not separate  $K \cap \Delta$
from $K' \cap \Delta$ in $\Delta$.

We stabilize $\s$ by tubing it along $\alpha$; denote the tube by $t$, the
surface obtained after tubing by $S(\s)$, and the complementary handlebodies
by $H_1'$ and $H_2'$ (with $K' \subset H_1'$).

Let $X'$ be the exterior of $K'$.  Since $K'$ is isotopic to $K$, $X' \cong
X$.  Note that $H_1'$ admits an obvious meridian disk that intersects $K'$
once (a component of $\Delta \cap H_1'$).  Note also that $K'$ is isotopic
into $S(\s)$ in $H_1'$.  Therefore, $S(\s)$ is a Heegaard surface for $X'$.
Since $g(S(\s)) = g(\s) + 1 = p+1$ and by assumption $p < g(X) = g(X')$, we
have that $S(\s)$ is a minimal genus Heegaard surface for $X'$.

We claim that $d(S(\s)) \leq 2$.  To prove this we will show that $H_1'$ and
$H_2'$ admit meridian disks that are disjoint from $K \subset S(\s)$.  In
$H_2'$ we take the compressing disk for the tube $t$.  For $H_1'$, let $D
\subset H_1$ be any meridian disk.  We will use $D$ to construct $D'$, a
meridian disk for $H_1'$, so that $D' \cap K = \emptyset$.  (Intuitively, we
construct $D'$ by pushing $D$ over $t$.)   Via isotopy we  may assume that $D$
intersects $N$ (if at all)  in disks $D_1,\dots,D_l$ (for some integer $l$)
that are parallel to $\Delta$,  and close enough to $\Delta$ so that $t$
intersects $D_i$ in the same pattern as it intersects  $\Delta$
($i=1,\dots,l$).  Note that $D$ cut open along $S(\s)$ has $2l+1$ components:
$l$ components inside $t$, $l$ components that intersect $K$, and  exactly one
other component, denoted $D'$.  It is easy to see that $D'$ is a meridian disk
for $H_1'$ disjoint from $K$.  Thus  $d(S(\s)) \leq 2$.

Let $\s'$ be any Heegaard surface for $X'$.  To estimate $d(\s')$ we apply
\cite[Corollary~3.5]{schtom} (with $S(\s)$ corresponding to $Q$ and $\s'$ to
$P$).  Then by \cite[Corollary~3.5]{schtom} one of the following holds:
\begin{enumerate}
\item $S(\s)$ is isotopic to a stabilization of $\s'$.  
\item $d(\s') \leq 2g(S(\s))$.
\end{enumerate}
In case~(1), since $S(\s)$ is a minimal genus Heegaard splitting, $S(\s)$ is
isotopic to $\s'$ (with no stabilizations).  Therefore $d(\s') = d(S(\s)) \leq
2 < 2(p+q)$.  On case~(2), $d(\s') \leq 2g(S(\s)) = 2(p+1) = 2(p+q)$.  As $X'
\cong X$, any Heegaard surface for $X$ has distance at most $2(p+q)$.
\end{proof}

\section{Calculating $g(X(m)^{(c)})$.}
\label{sec:genus}

Recall that we follow the notations in Definition~\ref{defns}.

\begin{lem}
\label{upperbound}
Let $K \subset M$ be a knot, $X$ the exterior of $K$ and $c \geq 0$ an
integer.  Denote $g(X)$  by $g$.

Then $g(X^{(c)}) \le g+c.$
\end{lem}

\begin{proof}
Note that $X^{(c)}$ is obtained from $X^{(c-1)}$ by drilling out a curve
parallel to $\partial X$.  Equivalently, we obtain $X^{(c-1)}$ by Dehn filling
a component of $\partial X^{(c)}$, and the core of the attached solid torus is
isotopic into $\partial X$.  This shows that the core of the solid torus is
isotopic to any  Heegaard surface of $X^{(c-1)}$, because one compression body
of the  Heegaard splitting is obtained from a regular neighborhood of
$(\partial X \cup  (\mbox{some components of } \partial X^{(c-1)} \setminus
\partial X)$  by adding some 1-handles.   In \cite{rieck}, this situation is
called a {\it good} Dehn filling, and it is shown that one of the following
holds:

\begin{enumerate}
\item $g(X^{(c)}) = g(X^{(c-1)})$. 
\item $g(X^{(c)}) = g(X^{(c-1)}) + 1$. 
\end{enumerate}

Hence we have 
$g(X^{(c)}) \le g(X^{(c-1)}) + 1$ in general. 
Since $g(X^{(0)}) = g(X) = g$, 
this implies $g(X^{(c)}) \le g+c$. 
\end{proof}

\begin{pro}
\label{pro:genus-x-super}
Let $M$ be a compact orientable manifold that does not admit a non-separating
surface, and $K \subset M$ a knot.  Let $c \geq 0$ be an integer. Denote $g(X)$
by $g$. Suppose that $X$ does not admit an essential surface $S$ with $\chi(S)
\geq  4 - 2(g+c)$, and that $K$ does not admit a $(g-1,c)$ position.  

Then $g(X^{(c)}) = g+c.$  
\end{pro}

\begin{proof}

The proof is an induction on $c$.  For $c = 0$ there is nothing to prove.
Fix $c>0$  as in the statement of the proposition and let $\s \subset X^{(c)}$
be a minimal genus Heegaard surface.  
Suppose that  $X$ does not admit an essential surface $S$ with 
$\chi(S) \geq  4 - 2(g+c)$, and that $K$ does not admit a 
$(g-1, c)$ position. 
By the inductive hypothesis we have 
$g(X^{(c-1)}) = g+c-1$.  
By the inequalities in the proof of Lemma~\ref{upperbound}, 
we have either 
$g(X^{(c)}) = g+c-1$, or 
$g(X^{(c)}) = g+c$.  
The proof is divided into the following two cases. 

\medskip

\noindent{\bf Case 1.}  $\s$ is strongly irreducible. 

\medskip

By Proposition~\ref{pro:bridge} one of the following holds:

\begin{enumerate}
\item $X$ admits an essential surface $S$ with $\chi(S) \geq 4 - 2g(X^{(c)})$.
\item $K$ admits a $(g(X^{(c)}) - b, b)$ position for some $b \geq c$.
\end{enumerate}

By Lemma~\ref{upperbound}, we have 
$4 - 2g(X^{(c)}) \ge 4-2(g+c)$. 
By assumption, $X$ does not admit an essential surface $S$ with $\chi(S) \geq
4 - 2(g+c)$, so~(1) above cannot happen and we may assume that $K$ admits a
$(g(X^{(c)}) - b, b)$ position for some $b \geq c$.  Since $b-c \geq  0$ we
can tube the Heegaard surface giving the $(g(X^{(c)}) - b, b)$ position
$b-c$ times to obtain a $(g(X^{(c)}) - b + (b-c), b-(b-c)) = (g(X^{(c)}) - c,
c)$ position.  Since $K$ does not admit a $(g-1, c)$ position, 
this implies that $g(X^{(c)}) - c \geq g$ and in particular $g(X^{(c)}) \ne
g+c-1$.  Hence we have $g(X^{(c)}) = g+c$. 

\medskip

\noindent{\bf Case 2.}  $\s$ is weakly reducible. 

\medskip 

By Casson and Gordon \cite{casson-gordon}, an appropriately chosen
weak reduction yields an essential surface $\widehat{F}$ (see
\cite[Theorem~1.1]{sedgwick} for a  relative version of Casson and Gordon's
Theorem).   Let $F$ be a connected component of $\widehat{F}$.  Since $F 
\subset X^{(c)} \subset M$ it separates and by \cite[Proposition~2.13]{cag}
$\s$ weakly reduces to $F$.  Note that $\chi(F) \geq \chi(\s) + 4$.

\medskip
\noindent{Claim.}  $F$ can be isotoped into $Q^{(c)}$ 
(recall the definition of $T$, $X'$ and $Q^{(c)}$ from the proof of
Proposition~\ref{pro:bridge}).

\begin{proof}[Proof of Claim.]
Assume for a contradiction this is not the case.  Since $F$ and $T$ are
essential, the intersection consists of a (possibly empty) collection of
curves that are essential in both surfaces.  Minimize $|F \cap T|$ subject to
this constraint.   If $F \cap X'$ compresses, then (since the curves of $F
\cap T$ are essential in $F$) so does $F$, contradiction.  Since $T$ is a
torus, boundary compression of $F \cap X'$ implies a compression (see, for
example, \cite[Lemma~2.7]{KRlocal}).  Finally, minimality of $|F \cap T|$
implies that no component of $F \cap X'$ is boundary parallel.  Thus, every
component of $F \cap X'$ is essential.  This includes the case $F \subset X'$
(in that case $F$ is essential in $X'$, else it would be parallel to $T$ and
isotopic into $Q^{(c)}$).  Since no component of $F \cap Q^{(c)}$ is a disk or
a sphere, $\chi(F \cap X')  \geq \chi(F) \geq \chi(\s) + 4$.   By
Lemma~\ref{upperbound}, we have  $\chi(\s) + 4 \geq 2 - 2(g+c) + 4 =6 -
2(g+c)$.  Hence $\chi(F \cap X') \ge 6 - 2(g+c)$, contradicting our assumption.
\end{proof} 

By \cite[VI.34]{jaco} $F$ is a vertical torus in $Q^{(c)}$.  First, if $F$ is
not parallel to a component of $\partial Q^{(c)}$ then  $F$ decomposes
$X^{(c)}$ as $X^{(p+1)} \cup_F D(c-p)$, where $p\geq 0$ is an integer and
$D(c-p)$ is a disk with $c-p$ holes cross $S^1$.  Note that since $F$ is not
parallel to a component of $\partial Q^{(c)}$, $c-p \geq 2$.  Therefore $p+1 <
c$ and by the inductive hypothesis  $g(X^{(p+1)}) = g+p+1$.  By Schultens
\cite{schultens}, $g(D(c-p)) = c-p$.

Next, suppose $F$ is boundary parallel in $Q^{(c)}$. 
Since $\s$ is minimal genus, $F$ is not parallel to a component of 
$\partial X^{(c)}$ \cite{sedgwick}. 
Hence $F$ is isotopic to $T$.  
This gives the decomposition $X^{(c)} = X' \cup_F
Q^{(c)}$.  Since $X' \cong X$, $g(X') =g$.  By \cite{schultens} $g(Q^{(c)}) =
c+1$.

Since $F$ was obtained by weakly reducing a minimal genus Heegaard surface,
\cite[Proposition~2.9]{cag}  (see also \cite[Remark~2.7]{schultens}) 
gives, in the first case:
 \begin{eqnarray*}
    g(X^{(c)}) &=& g(X^{(p+1)}) + g(D(c-p)) - g(F) \\
       &=& (g+p+1) + (c-p) - 1 \\ 
        &=&  g+c.
  \end{eqnarray*}
And in the second case:
 \begin{eqnarray*}
g(X^{(c)}) &=& g(X') + g(Q^{(c)}) - g(F) \\
&=& g + (c+1) - 1 \\
&=& g+c.
  \end{eqnarray*}
This completes the proof. 
\end{proof}

\begin{pro}
  \label{pro:genus}
Let $m \geq 1$ and $c \geq 0$ be integers and $\{K_i \subset M_i\}_{i=1}^m$
knots in closed orientable manifolds so that (for all $i$) $M_i$ does not
admit a non-separating surface.  Denote $E(K_i)$ by $X_i$ and $E(\#_{i=1}^m K_i)$ by $X$.
Let $g$ be an integer so that $g(X_i) \leq g$ for all $i$.

Suppose that no $X_i$ admits an essential surface $S$ with $\chi(S) \geq
\valchi$ and that no $K_i$ admit a $(g(X_i)-1,m+c-1)$ position.  Then we have:
$$g(X^{(c)}) = \s_{i=1}^m g(X_i)+c.$$
\end{pro}

\begin{rmkks}
  \begin{enumerate}
\item The proof for $m \geq 2$ is an induction of $(m,c)$ ordered
  lexicographically.  During the inductive step, $(m,c)$ is replaced by (say)
  $(m_1,c_1)$.  Since the complexity is reduced, $m_1 \leq m$.  However, $c_1 >
  c$. is possible.  We will see that if $c_1 >c$, then $c_1 = c+1$ and $m_1 <
  m$.  Thus $m_1 + c_1 \leq m+c$ and the condition `` no $X_i$ admits an
  essential surface $S$ with $\chi(S) \geq \valchi$'' holds when $m + c$ is
  replaced by $m_1 + c_1$.  The same holds for the condition ``no $K_i$ admit
  a $(g(X_i)-1,m+c-1)$ position''. 
  \item For $m \geq 2$, the proof is an application of the Swallow Follow Torus
  Theorem \cite[Theorem~4.1]{cag}.  In \cite[Remark~4.2]{cag} it was shown by
  means of a  counterexample that the Swallow Follow Torus Theorem does not
  apply  to $X^{(c)}$ when $m=1$.  Hence argument of
  Proposition~\ref{pro:genus} cannot be used to simplify the proof of
  Proposition~\ref{pro:genus-x-super}.
  \end{enumerate}
\end{rmkks}

\begin{proof}
The assumptions of Proposition~\ref{pro:genus-x-super} hold and so that
the proposition establishes the case $m=1$
(note that $4-2(1+c)g \le 4-2(c+g)$ holds). 
Hence we assume from now on that
$m\geq 2$. 

We induct on $(m,c)$ ordered lexicographically, where $m$ is the number of
summands and $c$ is the number of curves drilled.  Note that by 
Miyazaki \cite{miyazaki} $m$ is well defined (see \cite[Claim~1]{cag}).

By Lemma~\ref{upperbound}, and Inequality~(\ref{eq:upper-bound}) in section~1, 
we get: $g(X^{(c)}) \leq
g(X) + c = g(E(\#_{i=1}^m K_i)) + c \leq \s_{i=1}^m g(E(K_i)) + c \leq 
mg + c$.  Since $g \geq 2$ we have that $g(X(m)^{(c)}) \leq (m+c)g$.

By assumption, for all $i$, $X_i$ does not admit an essential surface
$S$ with $\chi(S) \geq 4-2(m+c)g$.  Hence by the Swallow Follow Torus
Theorem \cite[Theorem~4.1]{cag} any minimal genus Heegaard surface for
$X^{(c)}$ weakly reduces to a swallow follow torus $F$ giving the
decomposition $X = X_I^{(c_1)} \cup_F X_J^{(c_2)}$, with $I \subset
\{1,\dots,m\}$, $K_I = \#_{i \in I} K_i$, $K_J = \#_{i \not\in I} K_i$, $X_I =
E(K_I)$, $X_J = E(K_J)$, and $c_1 + c_2 = c+1$ (for details see the first
paragraph of Section~4 of \cite{cag}). 
Note that $I = \emptyset$ or $I = \{1,\dots,m\}$ are possible.  However,
at least one of $I,\ \{1,\dots,m\} \setminus I$ is not empty and by
symmetry we may assume $I \neq \emptyset$.

If $I = \{1,\dots,m\}$ then $c_2 \geq 2$.  Hence $c_1 < c$ and the inductive
hypotheses applies to $X_I^{(c_1)}\cong X^{(c_1)}$, showing that
$g(X_I^{(c_1)}) = \s_{i=1}^m g(X_i)+c_1$.  Since $X_J^{(c_2)}$ is homeomorphic
to a disk with $c_2$ holes cross $S^1$, by \cite{schultens} $g(X_J^{(c_2)}) =
c_2$. Amalgamation along $F$ gives (recall that $c_1 + c_2 = c+1$): 
\begin{eqnarray*}
  g(X^{(c)}) &=& g(X_I^{(c_1)}) + g(X_J^{(c_2)}) - g(F) \\
  &=& (\s_{i=1}^m g(X_i) + c_1) + c_2 - 1  \\
  &=&  \s_{i=1}^m g(X_i)+c.
\end{eqnarray*}

If $I \neq \{1,\dots,m\}$ then the number of summands in $K_I$ and $K_J$ are
$|I|$ and $m - |I|$ (respectively) and are both less than $m$.  By
construction $c_1 \leq c+1$ and $c_2 \leq c+1$ hence $m_1 + c_1 \leq m+c$ and
$m_2 + c_2 \leq 
m+c$.  By the inductive hypothesis
$g(X_I^{(c_1)}) = \s_{i \in I} g(X_i) + c_1$ and $g(X_J^{(c_2)}) =  \s_{i
  \not\in I} g(X_i)  + c_2$. Amalgamation along $F$ gives:
\begin{eqnarray*}
  g(X^{(c)}) &=& g(X_I^{(c_1)}) + g(X_J^{(c_2)}) - g(F) \\
  &=& ( \s_{i \in I} g(X_i) + c_1) + ( \s_{i \not\in I} g(X_i) + c_2) - 1  \\
  &=&  \s_{i=1}^m g(X_i)  + (c_1 + c_2) - 1 \\
  &=&  \s_{i=1}^m g(X_i)  + c.
\end{eqnarray*}
This proves the proposition in both cases.
\end{proof}

\section{Proof of Theorem~\ref{thm:additive}.}
\label{sec:proof}

Fix $g \geq 2$ and $n \geq 1$.  Let $\kk$ be the collection of all knots $K
\subset S^3$ so that:
\begin{enumerate}
\item $g(E(K)) \leq g$.
\item $X$ does not admits an essential surface $S$ with $\chi(S) \geq
  \valchin$.
\item $K$ does not admit a $(g(X) - 1, n)$ position.
\end{enumerate}

Applying Proposition~\ref{pro:genus} with $m \leq n$ and $c=0$, we see that
the knots in $\kk$ fulfill Condition~(2) of Theorem~\ref{thm:additive}.

Fix $h$, $2 \leq h \leq g$.   By \cite[Theorem~3.1]{mms} there exist
infinitely many knots $K$ with $g(X) = h$, admitting a Heegaard splitting of
distance greater than $\val$ (for $g=2$ this was  obtained independently by
Johnson  \cite[Lemma~4]{j}).  Let $K$ be such a knot and $X$ its exterior. By
Proposition~\ref{pro:essential}, since $X$ admits a Heegaard splitting with
distance greater than $2gn-2$, $X$ does not admits an essential surface $S$
with $\chi(S) \geq \valchin$.  By Proposition~\ref{pro:jt}, since $X$
admits a Heegaard splitting with distance greater than $2(h+n-1)$, $K$ does
not admit a $(g(X)-1, n) = (h-1,n)$ position.  We see that $K \in \kk$ and
hence, $\kk$ contains infinitely many knots $K$ with $g(X) = h$.  This proves
that $\kk$  fulfills condition~(1) as well.

Let $K \subset S^3$ be a knot with $g(E(K)) =h$.  As noted in the
introduction, any knot in $S^3$ is admissible (in the sense of \cite{crelle})
and therefore by~\cite[Theorem~1.2]{crelle}  there exists $N$ so that if $n
\geq N$ then $g(E(nK)) < ng(E(K))$.  This shows that $K \not\in \kk$ for $n
\geq N$.  Hence $K \not\in \cap_{i=1}^\infty \kk$.  As $K$  was arbitrary,
$\cap_{i=1}^\infty \kk = \emptyset$.

This completes the proof of Theorem~\ref{thm:additive}.

%% \bibliographystyle{plain}   %
%% \bibliography{general}   

\end{document}